\documentclass[leqno,final]{siamltex}
\usepackage{amsmath,amssymb,amsfonts}
\usepackage{graphicx}
\usepackage{epsfig}
\renewcommand{\a}{\alpha}
\renewcommand{\d}{\delta}
\newcommand{\e}{\varepsilon}
\newcommand{\W}{\Omega}
\renewcommand{\P}{\mathbb{P}}
\newcommand{\R}{\mathbb{R}}
\newcommand{\mJ}{\mathcal{J}}
\newcommand{\mJha}{\mJ_h^{\alpha}}
\newcommand{\mA}{\mathcal{A}}
\newcommand{\Th}{\mathcal{T}_h}
\renewcommand{\grad}{\nabla}
\newcommand{\ch}{\chi_{h}^{\alpha}}

\DeclareMathOperator*{\argmin}{argmin}
\DeclareMathOperator*{\sgn}{sgn}
\def \dx[#1]{\ensuremath{\operatorname{d}\!{#1}}}

\newtheorem{remark}{Remark}[section]

\begin{document}

\title{An enhanced finite element method for a class of variational problems exhibiting the Lavrentiev gap 
phenomenon\footnote{This work was partially supported by the NSF through grant DMS-1318486.}}

\author{Xiaobing Feng\thanks{Department of Mathematics, University of Tennessee, Knoxville, TN 37996, U.S.A.
(\tt{xfeng@math.utk.edu})}.
\and
Stefan Schnake\thanks{Department of Mathematics, University of Tennessee, Knoxville, TN 37996, U.S.A.
({\tt schnake@math.utk.edu})}.
}

\maketitle

\begin{abstract}
This paper develops an enhanced finite element method for approximating a class of variational 
problems which exhibit the \textit{Lavrentiev gap phenomenon} in the sense that the minimum 
values of the energy functional have a nontrivial gap when the functional is minimized on spaces $W^{1,1}$ 
and $W^{1,\infty}$. To remedy the standard finite element method, which fails to converge for 
such variational problems, a simple and effective cut-off procedure is utilized to design the 
(enhanced finite element) discrete energy functional. In essence the proposed discrete energy 
functional curbs the gap phenomenon by capping the derivatives of its input on a scale of $O(h^{-\a})$
(where $h$ denotes the mesh size) for some positive constant $\a$. A sufficient condition 
is proposed for determining the problem-dependent parameter $\a$. Extensive 1-D and 2-D
numerical experiment results are provided to show the convergence behavior and the 
performance of the proposed enhanced finite element method. 
\end{abstract}



\section{Introduction}\label{section1}
This paper concerns with finite element approximations of variational problems whose solutions (or minimizers) 
exhibit the so-called \textit{Lavrentiev gap phenomenon} - a defect from the singularities of the solutions.  
Such problems are often encountered in materials sciences, nonlinear elasticity, and image processing 
(cf. \cite{FHM2003,Winter1996,CO2009} and the references therein). These variational problems can be 
abstractly stated as follows: 
\begin{align} \label{variational_problem}
        u = \argmin_{v\in \mA} \mJ(v),
\end{align}
where the energy functional $\mJ:\mA\to\R\cup\{\pm\infty\}$ is defined by
\begin{align} \label{energy_functional}
        \mJ(v) = \int_{\W} f(\grad v,v,x) \dx[x].
\end{align}
Where $\W\subset \R^n$ is an open and bounded domain, $f:\R^n\times\R\times\W\to\R$, called the 
density function of $\mJ$, is assumed to be a continuous function, The space 
$\mA := \{v\in W^{1,1}(\W): v = g \text{ on } \partial\W  \}$ is known as the admissible set,  
$g\in L^1(\partial\W)$ is some given function.  

Let $\mA_\infty := \mA\cap W^{1,\infty}(\W)$. Since $\W$ is bounded, then $\mA \subset \mA_\infty$ 
and consequently there holds
\begin{align} \label{eq3}
   \inf_{v\in \mA_1} \mJ(v) \leq \inf_{v\in \mA_\infty}\mJ(v).
\end{align}
Problem \eqref{variational_problem} is said to exhibit {\em the Lavrentiev gap phenomenon} whenever
\begin{align}\label{gap_phenomenon}
        \inf_{v\in \mA_1} \mJ(v) < \inf_{v\in \mA_{\infty}}\mJ(v),
\end{align}
in other words, when the strict inequality holds in \eqref{eq3}.

The gap between the minimum values on both sides of \eqref{gap_phenomenon} suggests that the 
the minimizer of the left-hand side must have some singularity which causes the gap. Such a singularity 
often corresponds to a defect in a material or an edge in an image. It has been known in the literature 
\cite{FHM2003,Winter1996,CO2009} that the gap phenomenon could happen not only for nonconvex energy functionals 
but also for strictly convex and coercive energy functionals. As a result, it is a very complicated 
phenomenon to characterize and to analyze as well as to approximate (see below for details), because 
the gap phenomenon can be triggered by quite different mechanisms and the definition of {\em the 
Lavrentiev gap phenomenon} is a very broad concept which covers many different types of singularities. 
To the best of our knowledge, so far there is no known general sufficient conditions 
which guarantee the existence of the gap phenomenon. 

The simplest and best known example of the gap phenomenon is Mani\'a's 1-D problem \cite{Mania1934},
where one minimizes the functional
\begin{align}\label{eqn1:1}
\mJ(v) = \int_0^1 v'(x)^6\bigl(v(x)^3-x\bigr)^2 \dx[x]
\end{align}
over all functions $v\in W^{1,1}(0,1)$ satisfying $v(0)=0$ and $v(1)=1$.  By inspection it is easy 
to see that  $u(x)=x^{\frac13}$ minimizes \eqref{eqn1:1} with a minimum value zero.  However, 
it can be shown that the minimum over space $W^{1,\infty}(0,1)$ (i.e., the space of all Lipschitz functions)
is strictly larger than zero. As a result, Mani\'a's problem does exhibits {\em the Lavrentiev gap phenomenon}.
Notice that $u'(x)=\frac13 x^{-\frac23}$ which blows up rapidly as $x\to 0^+$. 
Moreover, a more striking property, which was proved by Ball and Knowles (cf. \cite{BK1987}), is
that if $u_j$ is a sequence of functions in $W^{1,q}(0,1)$ for $q\geq \frac32$ with $u_j(0) = 0$ and 
$u_j(1)=1$ such that $u_j\to u$ a.e. as $j\to\infty$, then $\mJ(u_j)\to\infty$ as $j\to\infty$. 
Since the finite element space $V^h_r$ (see section \ref{section2} for its definition) is a 
subspace of $W^{1,\infty}$, the above properties of the functional $\mJ$ 
imply that the standard finite element approximations to Mani\'a's problem 
must fail to approximate both the minimizer and the minimum value of the functional. Such a
conclusion was indeed verified numerically in \cite{FHM2003,Winter1996}, also see Figure \ref{fig2:1} 
for another numerical verification. 
This negative result immediately leads to the following two conclusions: first, variational problems 
which exhibit the gap phenomenon are difficult and delicate to approximate numerically; second, 
nonstandard numerical methods must be designed for such problems in order to correctly 
approximate both the minimizers and the minimum values. The failure of the standard finite 
element method suggests that in order to ensure the convergence of any numerical method which 
uses $V^h_r$ as the approximation space, one needs to construct a discrete energy functional $\mJ_h$ 
which necessarily does not coincide with $\mJ$ on the finite element space $V^h_r$.

As expected, there have been a few successful attempts to design 
convergent numerical methods for variational problems with the gap phenomenon.  
Below we only focus on discussing the methods which use conforming finite element methods to 
approximate variational problems with the Mani\'a-type gap phenomenon, by which we mean 
that the minimizers of the variational problems blow up in the $W^{1,\infty}$-norm.
but it is important to note that some gap phenomenon problems have been solved with the 
use of penalty and nonconforming finite element methods \cite{CO2009,CO2010,Ortner2011}.   

The first numerical method was proposed by Ball and Knowles in \cite{BK1987}. To handle 
the difficulty caused by the rapid blow-up in $W^{1,\infty}$-norm of the minimizer $u$,
they proposed to approximate $u$ and its derivative $u'$ simultaneously, an idea which 
is often seen in mixed finite element methods. Specifically, the authors proposed to minimize
the discrete energy functional
\begin{align} \label{BK}
\mJ_h^{BK}(v_h,w_h) = \int_{\W} f(w_h,v_h,x) \dx[x] 
\end{align}
under the constraint
\[
\|\phi(v_h'-w_h)\|_{L^1(\W)}\leq \e_h
\]
for some super-linear function $\phi$ over all functions $(v_h,w_h)\in V^h_1\times V^h_0$, 
where $V^h_0$ and $V^h_1$ denote respectively the discontinuous piecewise constant space 
and the continuous piecewise affine finite element space 
associated with a mesh $\Th$ of $\W$.
Where $\{\e_h\}$ is a sequence such that $\e_h \to 0$ as $h\to 0$. Notice that 
$\mJ_h^{BK}$ essentially has the same form as the original functional $\mJ$ after
setting $w_h=v_h'$. While this method works and is well-posed on the discrete level, 
the decoupling of $v_h$ and $v_h'$ adds an additional layer of unknowns which increases 
the complexity of the discrete minimization problem. Moreover, its generalization to higher
dimensions is not straightforward.
The other major numerical developments were carried out by Z. Li {\em et al.} in \cite{Li1992,Li1995,BL2006}.  
Their work has brought two similar methods: an element removal method and a truncation method. 
Here we only detail the truncation method and briefly mention the element removal method
because the latter is similar to the former and the truncation method is more closely related to our method to be 
introduced in this paper. Let $s\geq 1$ and $M_h > 0$. Define the discrete energy functional
\begin{align}\label{Li}
\mJ_h^{Li}(v_h) = \sum_{T\in\Th} \mJ_h^{Li}(v_h;T)
\end{align}
where
\begin{align*}
\mJ_h^{Li}(v_h;T) &= \min\left\{\mJ_h(v_h;T),\,M_h\left(1+\|\grad v_h \|_{L^s(T)}\right) \right\},\\
\mJ_h(v_h;T) &= \int_{T} f(\grad v_h,v_h,x) \dx[x].
\end{align*}
Here the truncation substitutes the contribution of $\mJ_h(v_h,T)$ by another constant if $v_h$ 
behaves ``poorly" on $T$. The element removal method simply discards (i.e., sets  $\mJ_h^{Li}(v_h,T)=0$ on)  
those ``bad" elements. Both methods are robust and calculate the minimum value 
of $\mJ$ over $\mA_\infty$ (assuming the minimizer $u$ uniquely exists). However, the determination of
$M_h$ and $s$ (or ``bad" elements) requires a litany of \textit{a priori} assumptions, some of which depend 
on the sought-after exact minimizer $u$. 

The goal of this paper is to introduce an effective and robust numerical method which slightly alters 
the standard finite element method by a novel and simple cut-off procedure. Our approach is motivated 
by the rationale that the standard finite element method fails to work because the magnitude of
the gradient $\nabla u_h$ becomes too large (independent of the magnitude of $u_h$, 
where $u_h$ stands for the standard finite element solution) near 
the singularity points. So the idea of our cut-off procedure is simply to limit the growth of $|\nabla u_h|$ 
to $\mathcal{O}(h^{-\alpha})$ order in the whole domain $\W$, the resulting discrete energy functional is then given by
\begin{equation}\label{discrete_energy}
\mJ_h^\alpha(w_h)=\int_\W f\bigl(\chi_h^\alpha(\nabla w_h), w_h,x\bigr)\, dx,
\end{equation}
where $\chi_h^\alpha(\cdot)$ denotes the cut-off function (see section \ref{section3} for its definition).
It is important to note that, unlike the truncation method of \cite{BL2006}, the choice of 
the crucial parameter $\alpha$ does not depend on any \textit{a priori} knowledge about the exact 
minimizer $u$, instead, it only depends on the structure of the energy density function $f$ and the space $\mA$. 
Moreover, we shall provide a sufficient condition, which is 
easy to use, for determining an upper bound for $\alpha$ to ensure the convergence.

The organization of the paper is as follows. In section \ref{section2} we introduce the notation, 
preliminary results such as finite element meshes and spaces. In section \ref{section3} we state the 
variational problems we aim to solve and the assumptions under which we develop our numerical method. 
We then define our finite element method with a help of the above cut-off procedure. We also present 
the alluded sufficient condition for determining an upper bound for $\alpha$ and demonstrate its
utility using Mani\'a's problem. In section \ref{section4} we provide some extensive numerical
experiment results for two specific application problems to gauge the performance of the proposed 
enhanced finite element method. Finally, the paper is ended with some concluding remarks in section 
\ref{section5}.

\section{Preliminaries}\label{section2}
Standard function and space notation will be adopted throughout this paper. For example, 
for an open and bounded domain $\W\subset\R^n$ with boundary $\partial\W$, 
let $W^{1,p}(\W)$ for $1\leq p\leq \infty$ denotes 
the Sobolev space consisting of functions whose up to first order weak derivatives are $L^\infty$-integrable 
over $\W$ and $\|\cdot\|_{W^{1,p}(\W)}$ denotes the standard norm on $W^{1,p}(\W)$.  
$\mA$ and $\mA_\infty$ have been introduced in section \ref{section1}, we also define the space
$\mA_p:=\mA\cap W^{1,p}(\W)$ for any $p\in (1,\infty)$.

Let $f$ and $\mJ$ be the same as in \eqref{energy_functional} and the variational problem to be considered
in the rest of this paper is given by \eqref{variational_problem}.
Suppose that problem \eqref{variational_problem} has a unique solution $u$ which exhibits 
{\em the Lavrentiev gap phenomenon} as defined in section \ref{section1}, that is, we assume 
there holds inequality 
\begin{align}\label{gap_phenomenon_1}
        \inf_{v\in \mA_1} \mJ(v) < \inf_{v\in \mA_{\infty}}\mJ(v).
\end{align}
So our primary goal is to construct an effective and robust finite element method to approximate $u$.

To this end, let $\Th$ be a quasi-uniform triangular (when $n=2$) or tetrahedral (when $n=3$) mesh of 
$\W$ with mesh parameter $h>0$. For a positive integer $r$, we define the finite element space $V^h_r$ 
on $\Th$ by
\begin{align}\label{eqn2:3}
	V^h_r := \bigl\{v_h\in \mA\cap C^0(\overline{\W}):\,  {v_h}|_ T \in\P_r(T) \quad \forall T\in \Th \bigr\},
\end{align}
where $\P_r(T)$ denotes the set of all polynomials whose degrees do not exceed $r$. 

\begin{figure}%
\centerline{
\includegraphics[width=3.2in,height=2.5in]{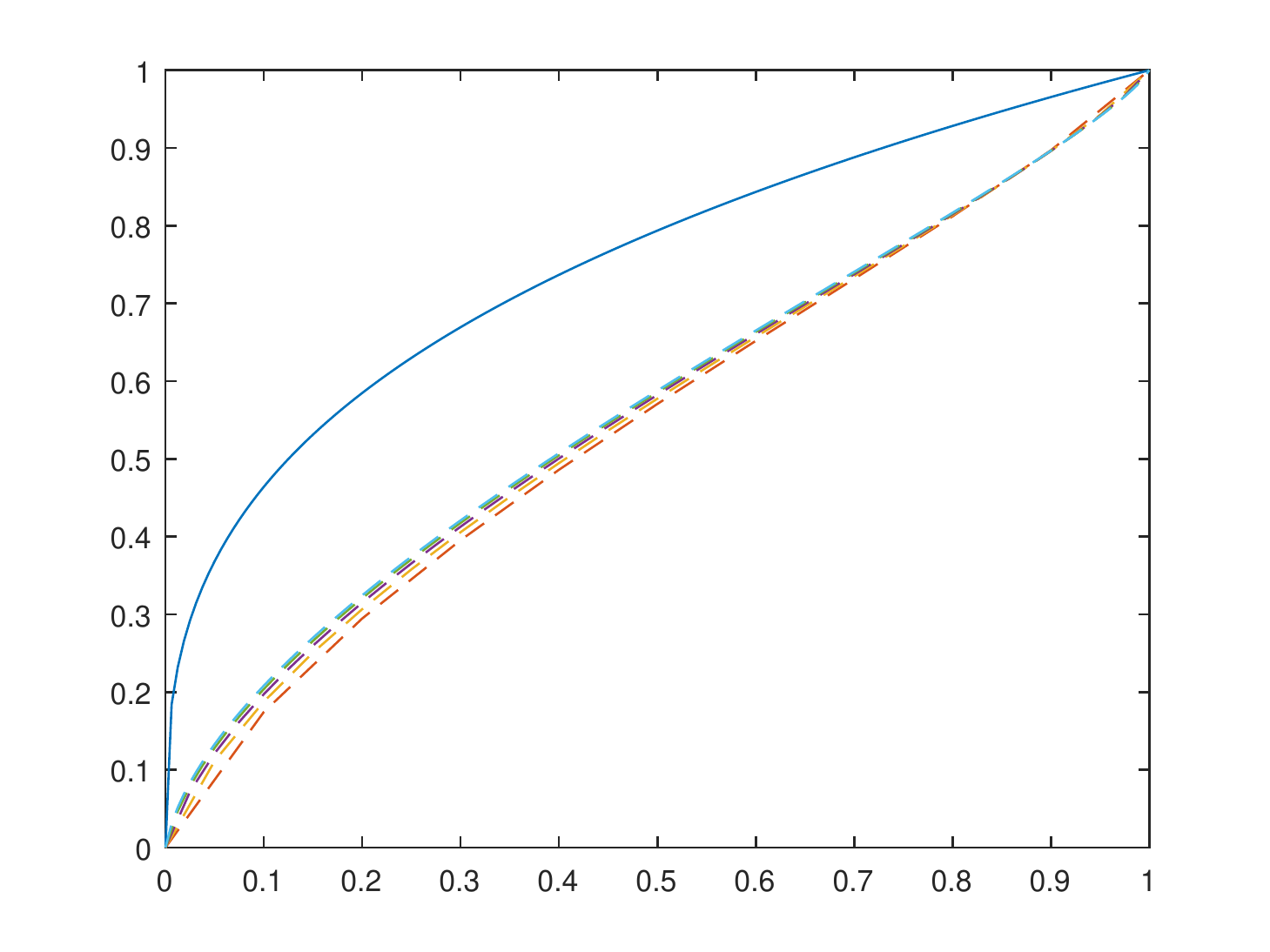}%
}
\caption{The standard finite element method applied to Mani\'{a}'s problem \eqref{eqn1:1}.  
The solid line is the true solution $u(x)=x^\frac13$ and the dashed lines are the finite element 
minimizers $u_h$ for $h=\frac{1}{N}$ where $N=10,20,40,80,160$. All minimizations were implemented 
by using the MATLAB minimization routine \texttt{fminunc} with initial function $u_0(x) = x$. }%
\label{fig2:1}%
\end{figure}

It is easy to see that $V^h_r\subset \mA_{\infty}\subset \mA$ for all $h>0$, then an obvious attempt 
to formulate a numerical method for variational problem \eqref{variational_problem} is the following 
standard finite element method which seeks $u_h\in V^h_r$ such that 
\begin{align}\label{FEM}
u_h = \argmin_{v_h\in V^h_r}\mJ(v_h).  
\end{align}
Unfortunately, this finite element method fails to give a convergent method because the method 
cannot give the true minimum value if \eqref{gap_phenomenon_1} holds and will not converge 
to the correct minimizer as the numerical test shows in Figure \ref{fig2:1}. 

To see the deeper reason, we note that for any $v\in\mA$ with $\mJ(v)<\infty$, the existence of a sequence 
of functions $\hat{v}_h\in V_h$ with $\hat{v}_h\to v$ in $\mA$ such that 
\begin{align}\label{eqn2:7}
\lim_{h\to 0}\mJ_h (\hat{v}_h) = \mJ(v)
\end{align}
is a key step to show convergence of the discrete minimizers. It is clear that \eqref{gap_phenomenon_1} 
implies that
\[
\mJ (\hat{u}_h) \geq  \inf_{v\in \mA_{\infty}}\mJ(v) > \mJ(u)
\]
for any $\hat{u}_h\in V^h_r$ with $\hat{u}_h\to u$ in $\mA$, which contradicts with \eqref{eqn2:7} 
for the minimizer $u$.  In fact, it was proved by C. Ortner in \cite{Ortner2011} that for a class of convex 
energies the convergence (to the exact solution) of the standard finite element method 
is equivalent to \eqref{gap_phenomenon_1} not holding (i.e., the gap phenomenon does not occur).  

\section{Formulation of the enhanced finite element method}\label{section3}

From the analysis given in the previous section we conclude that in order to construct
a convergent numerical method which uses $V^h_r$ as an approximation space, we must design a discrete 
energy functional $\mJ_h$ which should not coincide with $\mJ$ on the finite element space $V^h_r$.
In this section we shall construct a discrete energy functional $\mJ_h$ which meets this criterion
and provides a convergent (nonstandard) finite element method for problem \eqref{variational_problem}.   

Before introducing our method, let us give a heuristic discussion about why the gap phenomenon is appearing 
and how the existing methods assuage its effect.  Consider Mani\'{a}'s problem \eqref{eqn1:1}.  For 
any $v_h\in V_h$ (or in $\mA_{\infty}$) sufficiently approximating $u(x)=x^{\frac13}$, the quantity 
$(v_h^3-x)^2$ will be small but always nonzero. However, at the same time $|v_h'|$ will be very large 
near the origin.  If $|v_h'|$ is raised to a high enough power - six in this case - then the error 
of $(v_h^3-x)^2$ will be magnified to be so large that the quantity 
\[
\int_0^h (v_h')^6(v_h^3-x)^2 \dx[x]
\] 
will not vanish as $h\to 0$. For this reason, 
all of the existing methods were designed to dampen the effect of the 
derivative in the integral. The method of Ball and Knowles \cite{BK1987} weakly enforces $v_h'=w_h$ which 
allows the method to soften the effect of $v_h'$, where $v_h'$ has a singularity, and achieves 
convergence. The methods of Li {\em et al.} \cite{BL2006} leave the function $f$ unchanged, but remove 
or replace the functional value on the elements where something has gone wrong.  

With this in mind we now introduce a discrete energy functional which is much simpler and has a majority of the characteristics of the methods in \cite{Li1992,Li1995,BL2006}. Our approach is motivated by the 
belief that the standard finite element method fails to work because the magnitude of the gradient 
$\nabla u_h$ becomes 
too large (independent of the magnitude of $u_h$, where $u_h$ denotes the solution to \eqref{FEM}) 
near the singularity points. So our idea is simply to use a cut-off procedure to limit the growth 
of $|\nabla u_h|$ to $\mathcal{O}(h^{-\alpha})$ on the whole domain $\W$ in our discrete energy functional
$\mJ_h$. To this end, let $\a>0$, define the cut-off function $\ch:\R^n\to\R^n$ in the $i$th component by
\begin{align} \label{eqn3:3}
\left[\ch(s)\right]_i = \begin{cases} 
         s_i &\text{ if } |s_i| \leq h^{-\a} \\ 
         \mbox{sign}(s_i) h^{-\a} &\text{ if } |s_i| > h^{-\a}
                       \end{cases}, \qquad i=1,2,\ldots n.
\end{align}

It is clear that this function merely cuts the value of $s_i$ to a 
constant $\sgn(s_i) h^{-\alpha}$ if $|s_i|$ is too large. Then our discrete functional is simply
defined as
\begin{align}\label{eqn3:4}
\mJha(v_h) = \int_{\W} f\bigl(\ch(\grad v_h),v_h,x\bigr) \dx[x],
\end{align}
and our enhanced finite element method is defined by seeking $u_h\in V^h_r$ such that
\begin{align}\label{ehFEM}
u_h = \argmin_{v_h\in V^h_r}\mJha(v_h).  
\end{align}

\begin{remark}
Since our discrete energy functional $\mJha$ curbs the gap phenomenon by capping the derivative of 
its input on a scale of $\mathcal{O}(h^{-\alpha})$, spiritually it is similar to the truncation method 
of Li {\em et al.} \cite{BL2006}, but unlike the truncation method it keeps the dynamics 
of $f$ with respect to $v$ and $x$ much like Ball and Knowles' approach in \cite{BK1987}. 
Implementing the cut-off procedure is very simple and can be done by adding a few lines of code. Moreover,
unlike the truncation method, our enhanced finite element method does not require 
{\em a priori} knowledge about the exact minimizer $u$ of \eqref{variational_problem}.
Further adding to the simplicity is the existence of only one parameter $\a$ in the method. Here $\a$ 
controls the rate at which the cut-off grows and is the key for the convergence of the method.  
In general, $\a$ needs to be chosen in order to obtain equation \eqref{eqn2:7} for all $v\in\mA$ 
where $I_h v\in V^h_r$ is the finite element interpolant of $v$. Indeed, \eqref{eqn2:7} is the only
restriction we impose upon $\a$. A permissible range for $\a$, which guarantees convergence, depends on 
the structure of the density function $f$, so it is problem-dependent. Below we use 
Mani\'a's problem to demonstrate the process.
\end{remark}

We now derive an upper bound for $\alpha$ such that \eqref{eqn2:7} holds for $v\in\mA$. 
For a fixed $v\in\mA$, notice that $\mJ(v)<\infty$, let
$I_h v\in V^h_r$ denote the finite element interpolant of $v$. We want to find an uppper bound 
for $\alpha$ such that $\mJha(I_h v) \to \mJ(v)$ as $h\to 0$ because this will guarantee \eqref{eqn2:7} for $v$. 
Let $\d>0$. Adding and subtracting $v$, using Young's inequality with weight $h^{\d}$, 
and using the definition of $\ch$ we get 
\begin{align*}
&\mJ_h^{\alpha}(v_h) = \int_0^1 (\ch(v_h'))^6(v_h^3-x)^2 \dx[x] \\
&\quad= \int_0^1 (\ch(v_h'))^6(v_h^3-v^3+v^3-x)^2 \dx[x] \\ 
&\quad\leq \int_0^1 (1+h^{-\d})(\ch(v_h'))^6(v_h^3-v^3)^2 \dx[x] + \int_0^1 (1+h^{\d})(\ch(v_h'))^6(v^3-x)^2\dx[x] \\
&\quad\leq \int_0^1 (1+h^{-\d})h^{-6\a}(v_h^3-v^3)^2 \dx[x] + \int_0^1 (1+h^{\d})(v_h')^6(v^3-x)^2\dx[x] \\
&\quad=: A_1^h + A_2^h.
\end{align*}
Since $v^3-x$ factor now has no error, multiplying by $(v'_h)^6$ does not have a magnification effect 
which is the source of the gap phenomenon, it can be shown that \cite{FS2016} $A_2^h\to \mJ(v)$ 
as $h\to 0$, then we have \eqref{eqn2:7} provided $A_1^h$ vanishes.  We claim that $A_1^h$ vanishes 
as $h\to 0$ for $0<\a<\frac16$. The proof of the assertion goes as follows.
By H\"{o}lder's Inequality we have
\begin{align*}
A_1^h &= (1+h^{-\d})h^{-6\a}\int_0^1 (v_h^3-v^3)^2 \dx[x] \\
& = (1+h^{-\d})h^{-6\a}\int_0^1 (v_h-v)^2(v_h^2+v_h v + v^2)^2 \dx[x] \\
&\leq (1+h^{-\d})h^{-6\a}\|v_h-v\|_{L^2(\W)}^2\|(v_h^2+v_h v + v^2)^2\|_{L^{\infty}(\W)}.
\end{align*}
Since $v_h=I_h v$ we have that $v_h$ is uniformly bounded in $h$ and 
$\|v_h-v\|_{L^2(\W)}^2  = \mathcal{O}(h)$.  Thus
\begin{align*}
0\leq A_1^h \leq \|(v_h^2+v_h v + v^2)^2\|_{L^{\infty}(\W)}(1+h^{-\d})h^{1-6\a}
\end{align*}
Since $\a < \frac16$ we may choose $\d<1-6\a$ such that $A_1^h\to 0$ as $h\to 0$ and 
we have \eqref{eqn2:7}. Clearly,  this range of $\a$ does not depend on the solution $u$ 
but only on the form of $f$ and the regularity of the space $\mA$. We regard this property 
as one crucial advantage of our method. 

\section{Numerical experiments}\label{section4}

In this section we present some numerical experiment results for two variational
problems which are known to exhibit the gap phenomenon. The first problem is 
Mani\'{a}'s 1-D problem which has been seen in the previous sections; 
the second problem, which was proposed by Foss in \cite{Foss2003}, is a 2-D variational 
problem from nonlinear elasticity. For each of the two test problems we solve it by using our 
enhanced finite element method with linear element (i.e., $r=1$), and we solve 
the minimization problem \eqref{ehFEM} by using the MATLAB minimization function \texttt{fminunc}.  
We first demonstrate the convergence of the numerical method, we then numerically evaluate the effect 
and sharpness of the parameter $\a$, and compare with the standard finite element method 
(which is known to be divergent).  We also numerically compute the rate of convergence for 
$u-u_h$ although no theoretical rate convergence has yet been proved for the numerical method.

\subsection{Mani\'{a}'s 1-D problem}\label{subsection4.1}

Once again, the energy functional of Mani\'{a}'s 1-D problem is given by \eqref{eqn1:1}.  
A uniform mesh $\Th$ with mesh size $h$ and the linear finite element are used in the test. 
As mentioned above, we solve the resulting minimization problem \eqref{ehFEM} 
by using the MATLAB minimization function \texttt{fminunc} with initial function 
$u_0(x)=x$.

\begin{figure}%
\centerline{
\includegraphics[width=3.2in,height=2.4in]{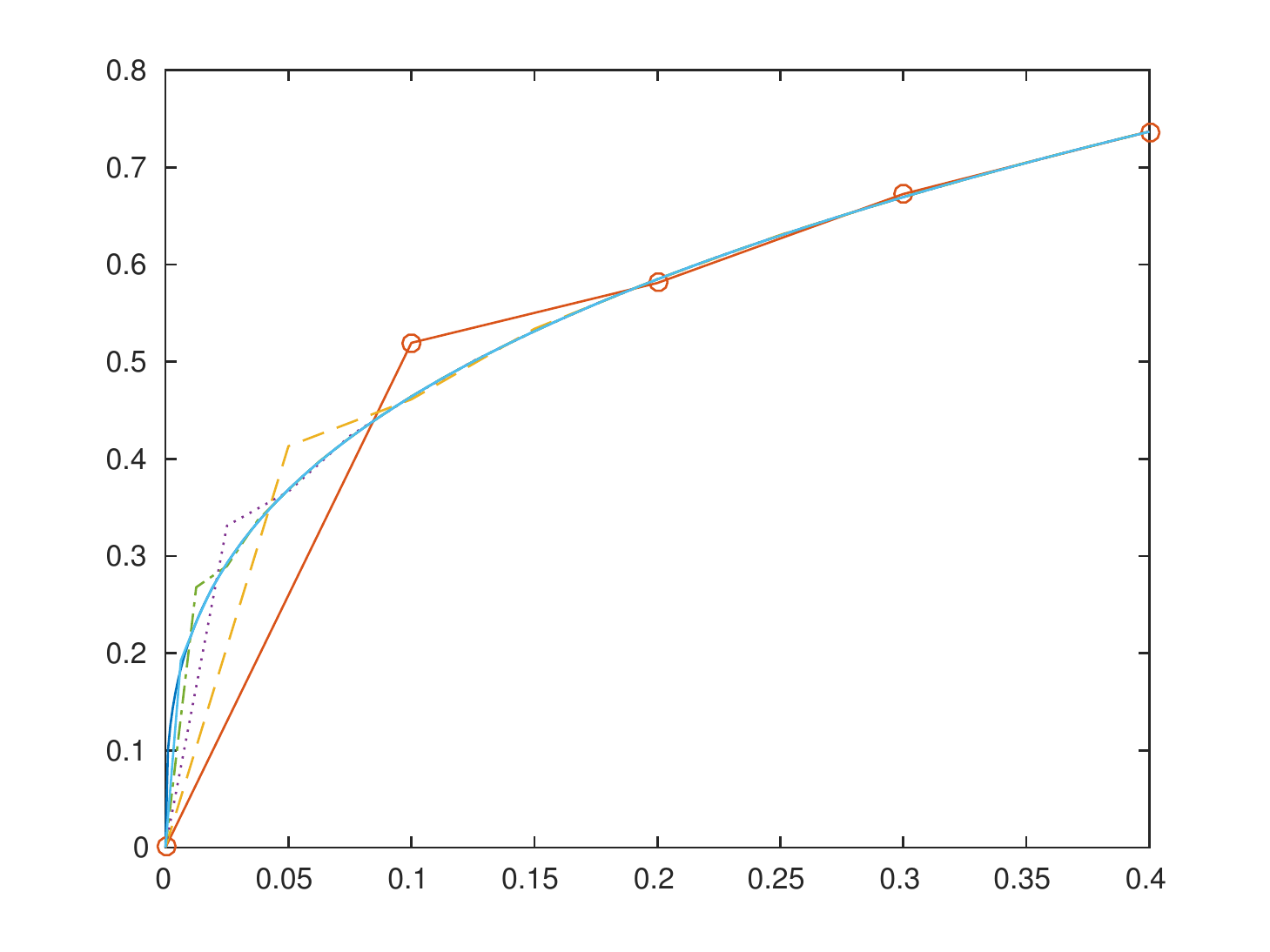}%
}
\caption{The graphs of the computed minimizers/solution $u_h$ of the enhanced FEM applied to Mani\'{a}'s 
problem \eqref{eqn1:1} with parameter $\a=\frac14$ from $x=0$ to $x=0.4$.  The solid line is the exact
solution $u(x)=x^{\frac13}$ and the dashed and circled lines are the minimizers $u_h$ for $h=\frac{1}{N}$ 
where $N=10,20,40,80,160$.  All minimizations were implemented by using the MATLAB minimization 
function \texttt{fminunc} with initial function $u_0(x) = x$.}%
\label{fig5:1}%
\end{figure}

\begin{figure}%
\centerline{
\includegraphics[width=3.2in,height=2.4in]{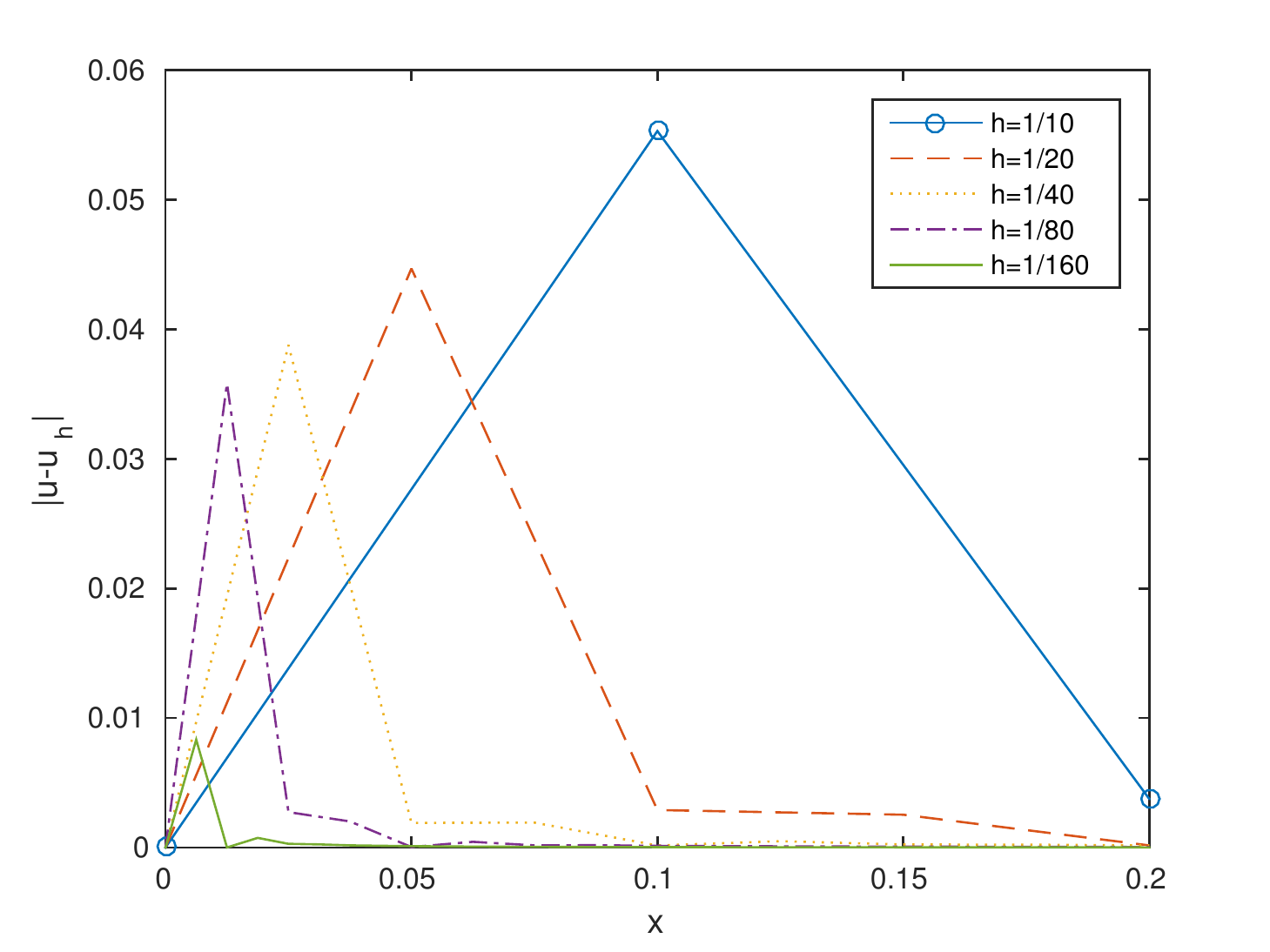}%
}
\caption{The graphs of the errors $|u_h-u|$ of the enhanced FEM applied to Mani\'{a}'s 
problem \eqref{eqn1:1} with parameter $\a=\frac14$ from $x=0$ to $x=0.2$ for $h=\frac{1}{N}$ 
where $N=10,20,40,80,160$.  All minimizations were implemented by using the MATLAB minimization 
function \texttt{fminunc} with initial function $u_0(x) = x$.}%
\label{fig5:3}%
\end{figure}

Figure \ref{fig5:1} displays the computed solutions (minimizers) $u_h$ with various mesh size $h$ 
along with the exact solution $u(x) = x^{\frac13}$. The parameter $\a=\frac14$ is used for the tests.  
It is clear that the solutions $u_h$ are correctly approximating $u$.  Figure \ref{fig5:3} shows 
the behavior of the absolute value of the error function $u-u_h$. As expected, we see that
the location where the biggest error occurs moves closer to the singularity point $x=0$ 
of $u$ as the mesh size $h$ gets smaller.

\begin{figure}%
\centerline{
\includegraphics[width=3.2in,height=2.4in]{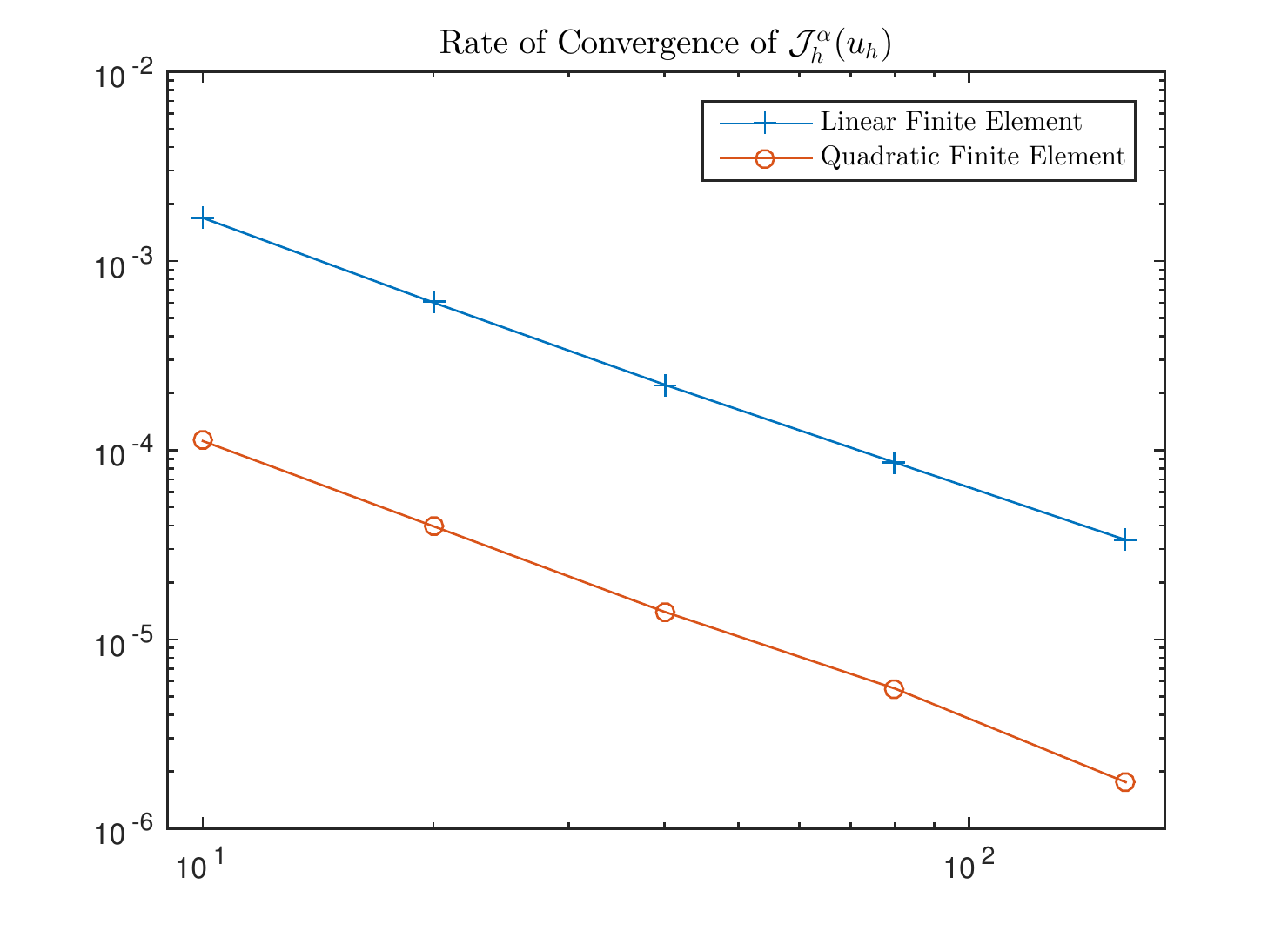}%
}
\caption{The rate of convergence of $\mJha(u_h)$ where $u_h$ is the solution to enhanced FEM applied to Mani\'{a}'s 
problem \eqref{eqn1:1} with parameter $\a=\frac14$ for $h=\frac{1}{N}$ 
where $N=10,20,40,80,160$. Plotted are the rates for the linear and quadratic finite element spaces. All minimizations were implemented by using the MATLAB minimization 
function \texttt{fminunc} with initial function $u_0(x) = x^{1/2}$.}%
\label{fig5:10}%
\end{figure}

For a more detailed look, we also record the $L^{\infty}$-norms of the error $u-u_h$ and compute the 
rate of convergence in Table \ref{table5:1}. Clearly, the table shows the convergence of the computed 
solutions $u_h$. As a comparison and to see that these approximations would not be found using the 
standard finite element method, a comparison of the values of $\mJ$ and $\mJha$ at $u_h$, $I^1_hu$, and $I_h^2(u)$
is given in Table \ref{table5:2}, where $I_h^1$ and $I_h^2$ are the piecewise linear and quadratic interpolants respectively.  
 We see here that $\mJha$ correctly captures the dynamics needed 
to obtain a convergent sequence of solutions $u_h$ while the sequences $\{\mJ(u_h)\}$ 
and $\{\mJ(I_h u)\}$ do not.  In addition $\{\mJha(I_h^1 u)\}$ and $\{\mJha(I_h^2 u)\}$ 
converge with the same rate, $\mathcal{O}(h^{1.5})$.  Thus employing higher order elements 
on this problem will not result in a larger convergence rate.  To make this clear we plot the convergence rate of the numerical minimizers $u_h$ of $\mJha$ for linear and quadratic elements in Figure \ref{fig5:10}. Note both elements observe the same convergence rate of 1.5.

\begin{table}[h]%
\begin{center}
\begin{tabular}{| c | c | c | c | c | c | }\hline
$ h $ & 1/10 & 1/20 & 1/40 & 1/80 & 1/160  \\ \hline
$\| u-u_h \|_{L^{\infty}}$ & 5.53e-2 & 4.50e-2 & 3.88e-2 & 3.59e-2 & 8.32e-3  \\ \hline
rate & - & 0.30 & 0.20 & 0.11 & 2.10  \\ \hline
\end{tabular}
\caption{The $L^{\infty}$ errors between $u$ and $u_h$ where $u_h$ are the solutions of 
the enhanced FEM applied to Mani\'{a}'s problem \eqref{eqn1:1} with parameter $\a=\frac14$.}
\label{table5:1}%
\end{center}
\end{table}

\begin{table}[h]%
\begin{center}
\begin{tabular}{| c | c | c | c | c | c | }\hline
$ h $ & 1/10 & 1/20 & 1/40 & 1/80 & 1/160  \\ \hline
$\mJ(u_h)$   & 8.23e-1 & 1.64 & 3.28 & 6.56 & 13.1 \\ \hline
$\mJha(u_h)$ & 1.68e-3 & 6.02e-4 & 2.22e-5 & 8.59e-6 & 3.31e-5   \\ \hline
$\mJ(I_h^1 u )$  & 7.19e-1 & 1.52 & 3.04 & 6.09 & 12.9   \\ \hline
$\mJha(I_h^1 u)$ & 2.41e-3 & 8.63e-4 & 3.09e-4 & 1.10e-4 & 3.91e-5   \\ \hline
$\mJ(I_h^2 u )$  & 1.16 & 2.31 & 4.62 & 9.24 & 18.5   \\ \hline
$\mJha(I_h^2 u)$ & 1.71e-4 & 6.03e-5 & 3.09e-4 & 7.54e-6 & 2.67e-6   \\ \hline
\end{tabular}
\caption{The functional values $\mJ$ and $\mJha$ for $u_h$, $I_h^1 u$, and $I_h^2 u$, where $u_h$ are the solutions 
of the enhanced FEM applied to Mani\'{a}'s problem \eqref{eqn1:1} with parameter $\a=\frac14$, 
and $I_h^1 u$ and $I_h^2 u$ is the piecewise linear and quadratic nodal interpolant of the exact solution/minimizer $u$.}
\label{table5:2}%
\end{center}
\end{table}

\begin{figure}[h]%
\begin{center}
\includegraphics[width=3.0in,height=2.5in]{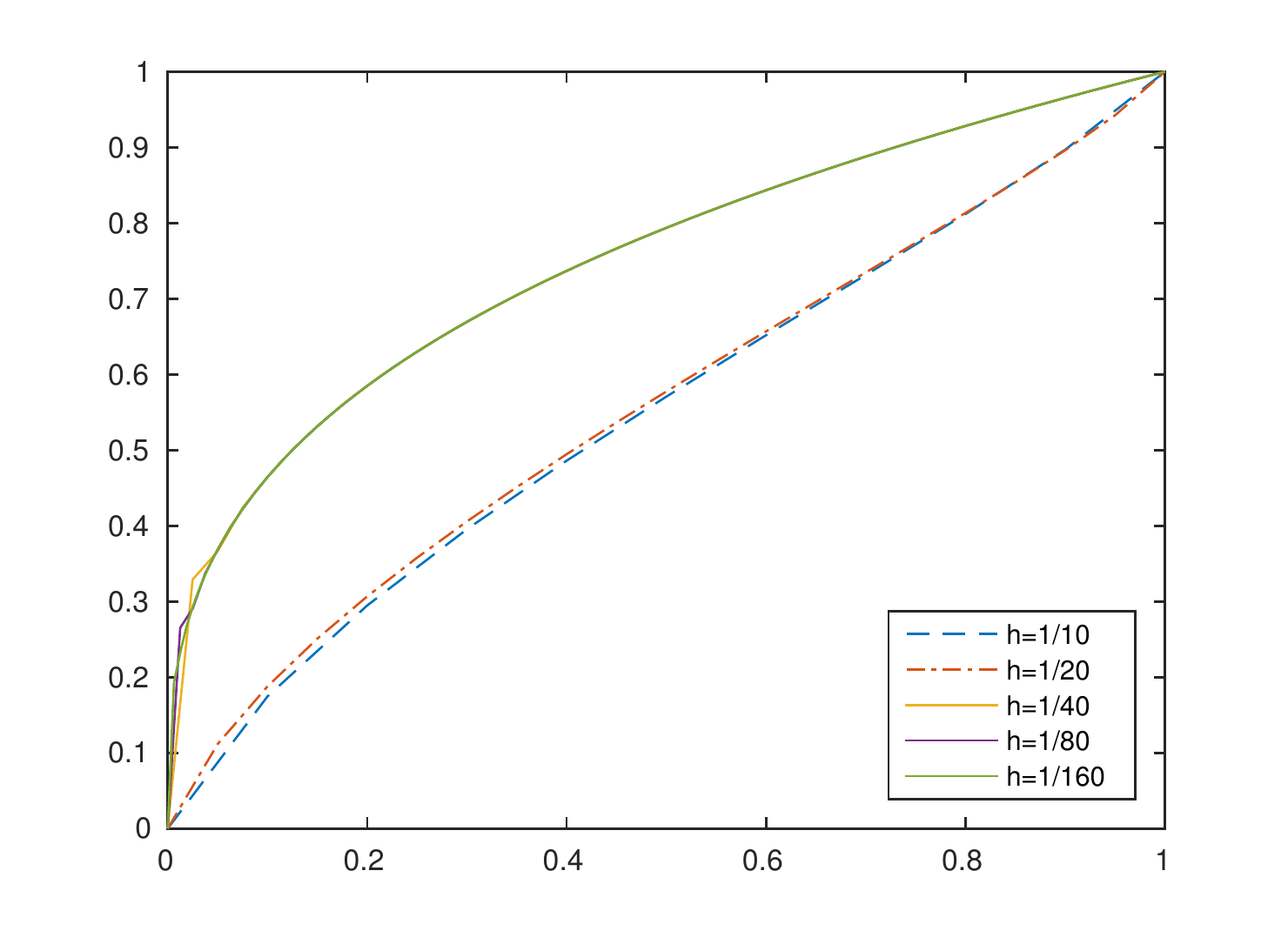}%
\caption{The graphs of the computed solutions/minimizers $u_h$ of the enhanced FEM applied to the 
Mani\'{a}'s problem \eqref{eqn1:1} with parameter $\a=\frac27$ for $h=\frac{1}{N}$ where $N=10,20,40,80,160$. 
The dotted lines are for $N=10$ and $20$ while the solid lines are for $N=40,80,$ and $160$. 
All minimizations were implemented by using the MATLAB minimization function \texttt{fminunc} 
with initial function $u_0(x) = x$.}%
\label{fig5:4}%
\end{center}
\end{figure}

Finally, we examine the role of the parameter $\a$. In section \ref{section3} we show that 
$\a < \frac12$ is sufficient to ensure \eqref{gap_phenomenon_1} for all $v\in\mA$ with finite 
energy.   Our numerical tests show that for any $\a<1/2$ the enhanced
finite element method converges for Mani\'{a}'s problem, and the convergence of $|\mJha(u_h)-\mJ(u)|\to 0$ 
diminishes as $\a\to \frac12$. So $\a^*:=\frac12$ seems a critical point for the choice of $\a$ for 
linear, quadratic, and higher order nodal finite elements.  
It must be noted that taking $\a$ close to $\a^*$ is not a good idea. Notice that the Euler-Lagrange 
equation of \eqref{eqn1:1} is a nonlinear equation. To solve the nonlinear equation, a mesh restriction 
$h<h'$ is expected and it takes up the most of the total CPU time for solving the nonlinear equation. 
This mesh restriction is expected to depend on $\a$.  To see this, let 
\[
\tilde{u}_h = \argmin_{v_h\in V_h} \mJ(u_h)
\] 
be the solution to the standard finite element method. Suppose that $\a$ is close to $\frac12$,
we observe that $\mJha(\tilde{u}_h)\approx\mJ(\tilde{u}_h)$. While $\mJha(I_h u)$ indeed converges 
to $\mJ(u)$ the convergence is very slow. Since the upper bound $h'$ must be chosen such that
for all $h<h'$ we have
\[
\mJha(I_h u) < \mJha(\tilde{u}_h), 
\]
so $h'$ must be extremely small and approaches $0$ as $\a\to \frac12$.  On noting the fact that
for all $h\geq h'$ a small perturbation of $\tilde{u}_h$ will be a minimizer of $\mJha$ over $V_h$, 
 we see that $\a$ must be chosen carefully in order to guarantee that we can obtain good 
numerical solutions with any mesh sizes $h<h'$. To show this important detail graphically, 
Figure \ref{fig5:4} displays the computed solutions/minimizers $u_h$ to $\mJha$ with $\a=\frac27$. 
We observe that for $h=\frac{1}{10}$ and $h=\frac{1}{20}$, $u_h$ do not approximate $u$ well, 
but for $h=\frac{1}{40}$, $h=\frac{1}{180}$ and $h=\frac{1}{160}$, $u_h$ gives much more accurate 
approximations.

\subsection{Foss' 2-D problem}\label{subsection4.2}

We now consider a 2-D variational problem which exhibits the Lavrentiev gap phenomenon. It arises from 
nonlinear elasticity and was first studied by M. Foss in \cite{Foss2003}, and its numerical approximation 
was investigated by Li {\em et al.} in \cite{BL2006}.

Let $\W=(0,1)\times(\frac{3}{2},\frac{5}{2})$, the energy functional of Foss' problem is given by  
\begin{align}\label{eqn5:3}
\mJ(v) = 66\Bigl(\frac{13}{14} \Bigr)^{14} \int_{\W} \Bigl(\frac{y}{y-1}\Bigr)^{14}
|u|^{\frac{14-3y}{y-1}}\bigl(|u|^{\frac{y}{y-1}}-x\bigr)^2(u_x)^{14} \dx[x]\dx[y],
\end{align} 
and the admissible set is $\mA = \{u\in W^{1,1}(\W) : u(0,\cdot) = 0 \text{ and } u(1,\cdot) = 1\}$.  
It was shown by Foss \cite{Foss2003} that
\[
0=\inf_{v\in\mA}\mJ(v)  <  \inf_{v\in\mA_{\infty}}\mJ(v)=1,
\]
which proves that $\mJ$ does exhibit the gap phenomenon. Moreover, the minimizer of $\mJ$ over 
$\mA$ is given by $u(x,y) = x^{\frac{y-1}{y}}$, but the problem does not attain its minimum value 
in $\mA_{\infty}$.  

We apply our enhanced finite element method with $\a=\frac16$ to solve Foss' problem. In order to 
generate a reasonably good initial guess for using the MATLAB minimization function \texttt{fminunc},
we first compute 
\begin{align} \label{eqn5:4}
	\tilde{u}_h = \argmin_{v_h\in V_h} \mJ(v_h)
\end{align}
using the MATLAB routine \texttt{fminunc} with initial guess $u(x,y) = x$ and then use 
$\tilde{u}_h$ as an initial condition for solving 
\begin{align} \label{eqn5:5}
	u_h = \argmin_{v_h\in V_h}\mJha(v_h)
\end{align}
using the same MATLAB routine \texttt{fminunc}. 

\begin{figure}%
\includegraphics[width=0.5\columnwidth]{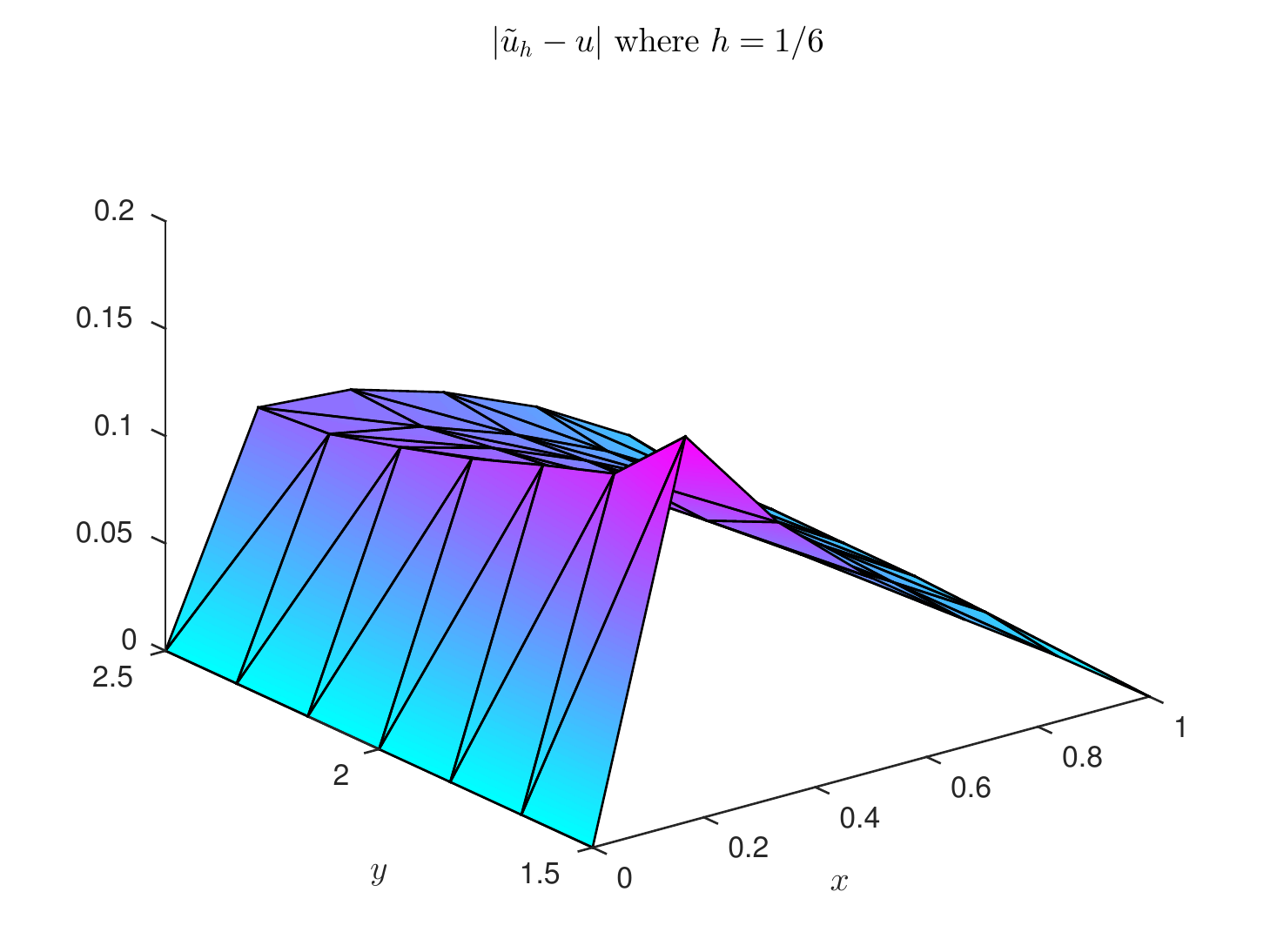}
\includegraphics[width=0.5\columnwidth]{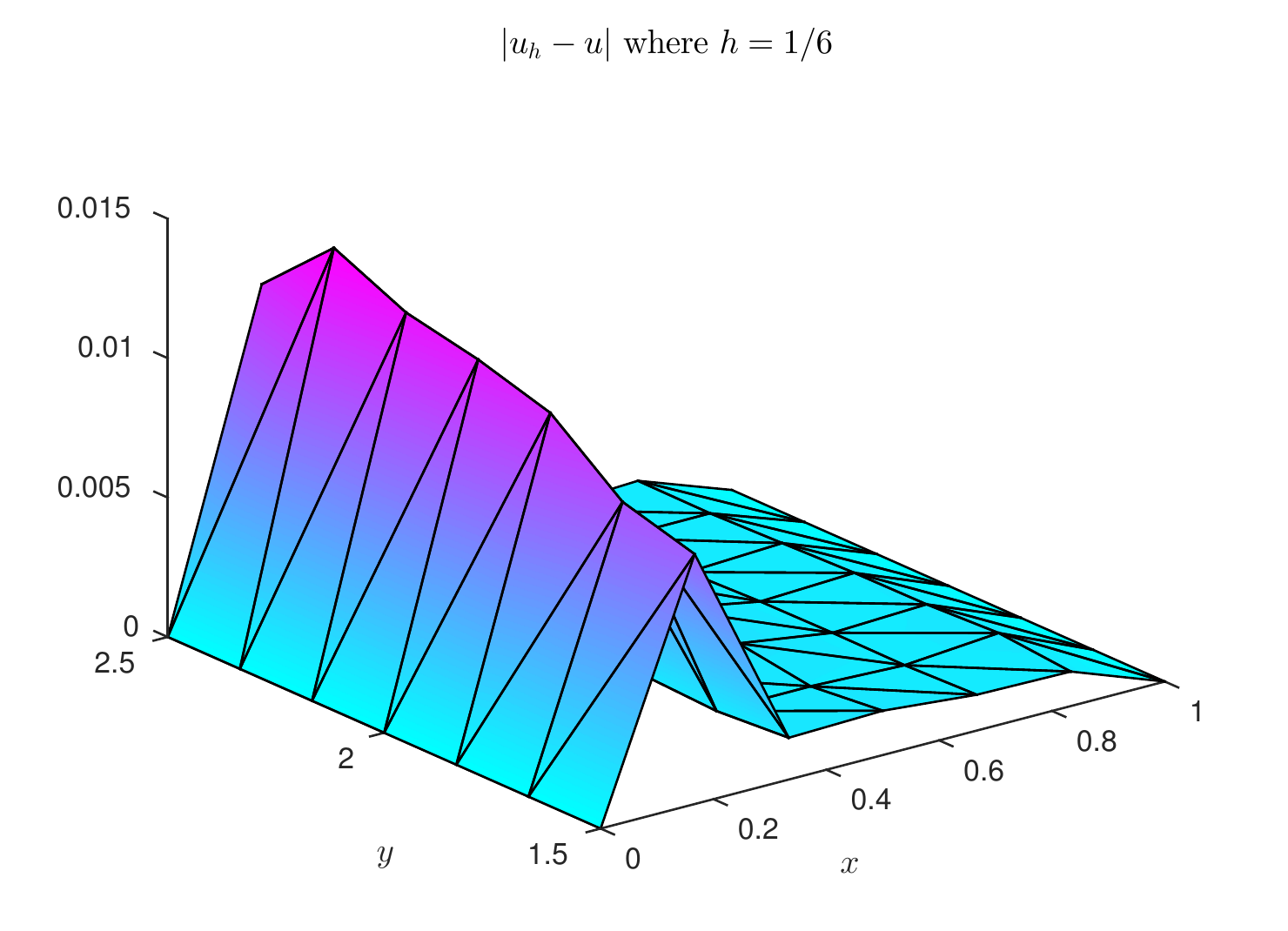}

\includegraphics[width=0.5\columnwidth]{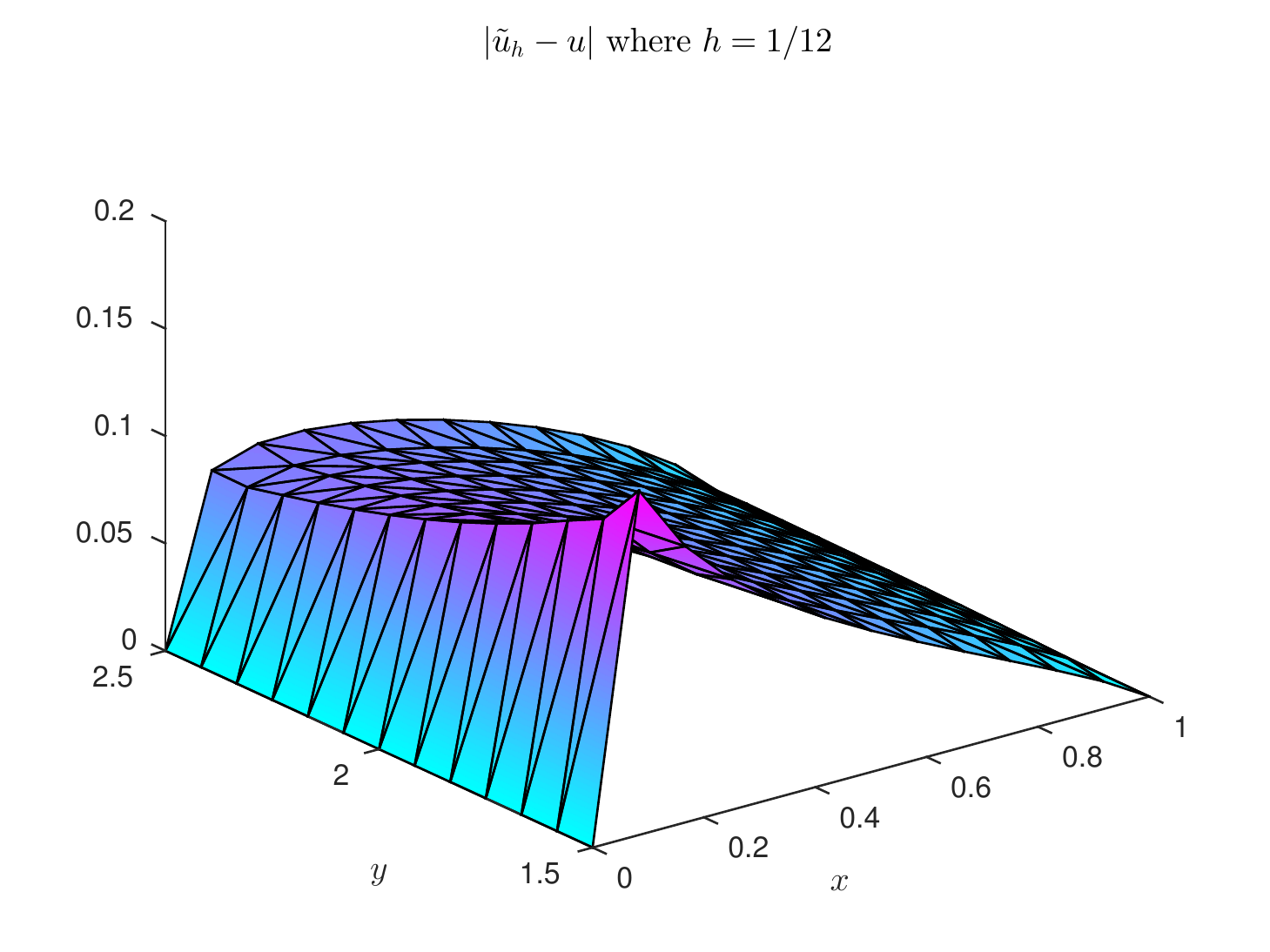}
\includegraphics[width=0.5\columnwidth]{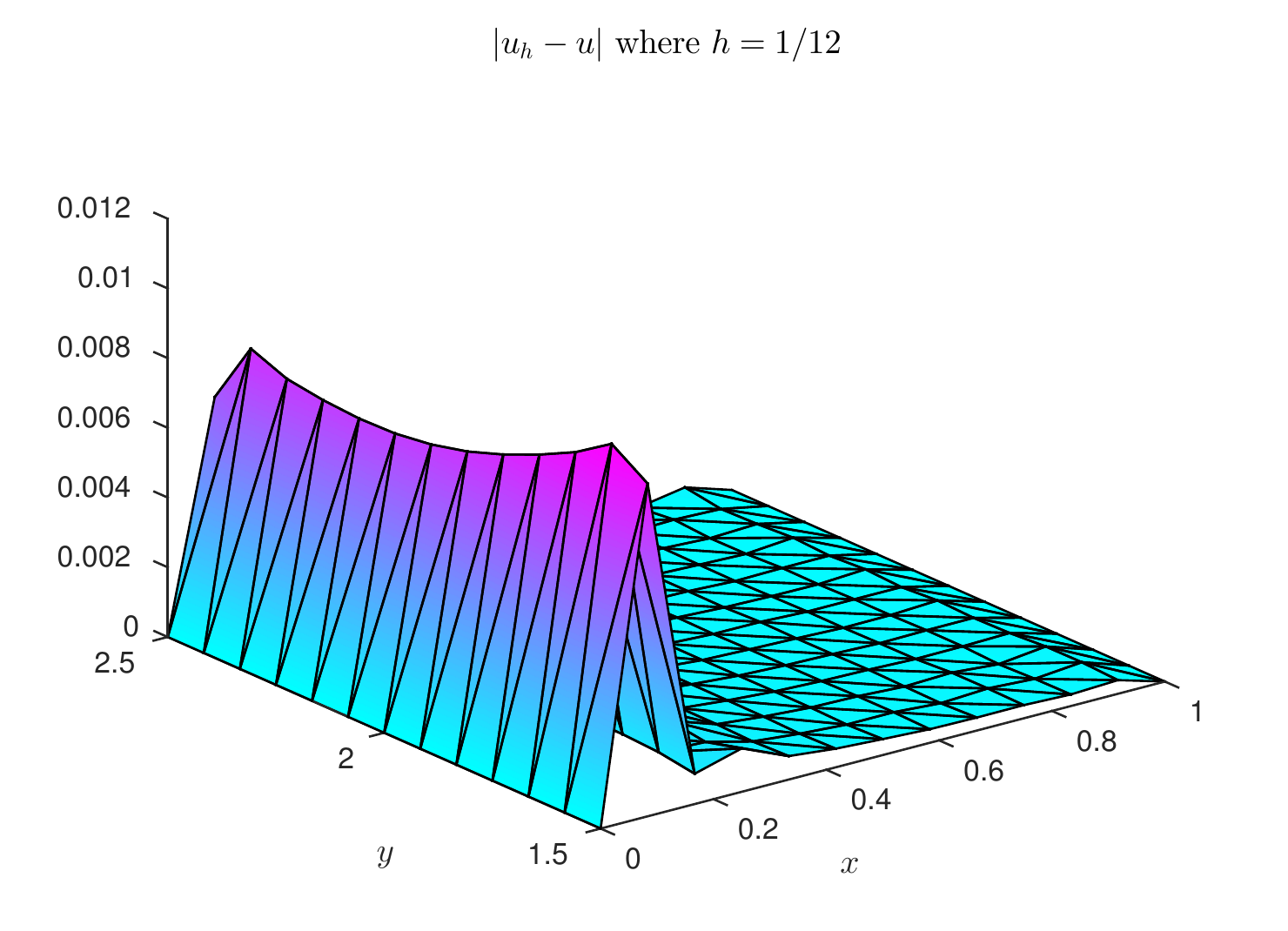}

\includegraphics[width=0.5\columnwidth]{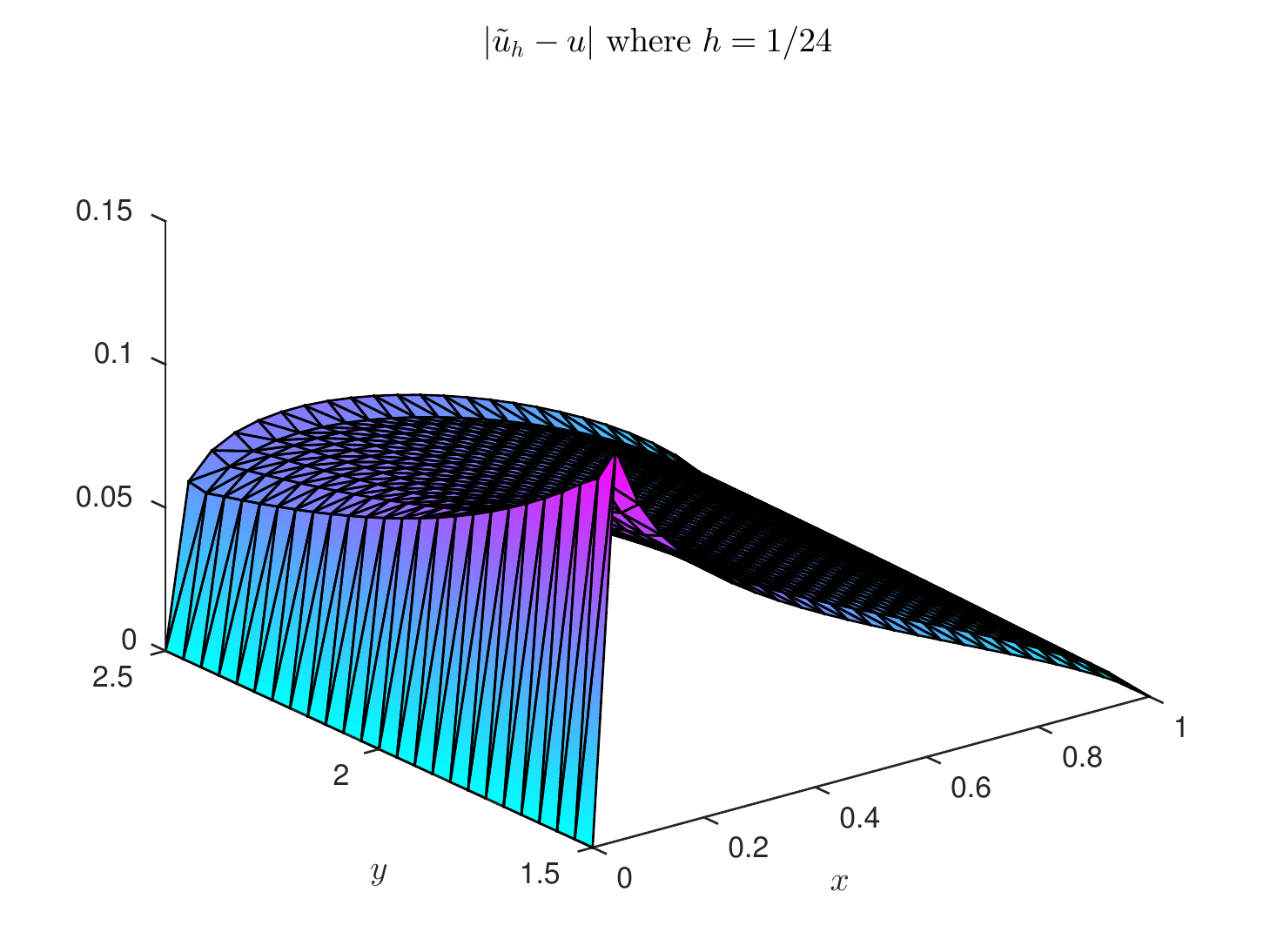}
\includegraphics[width=0.5\columnwidth]{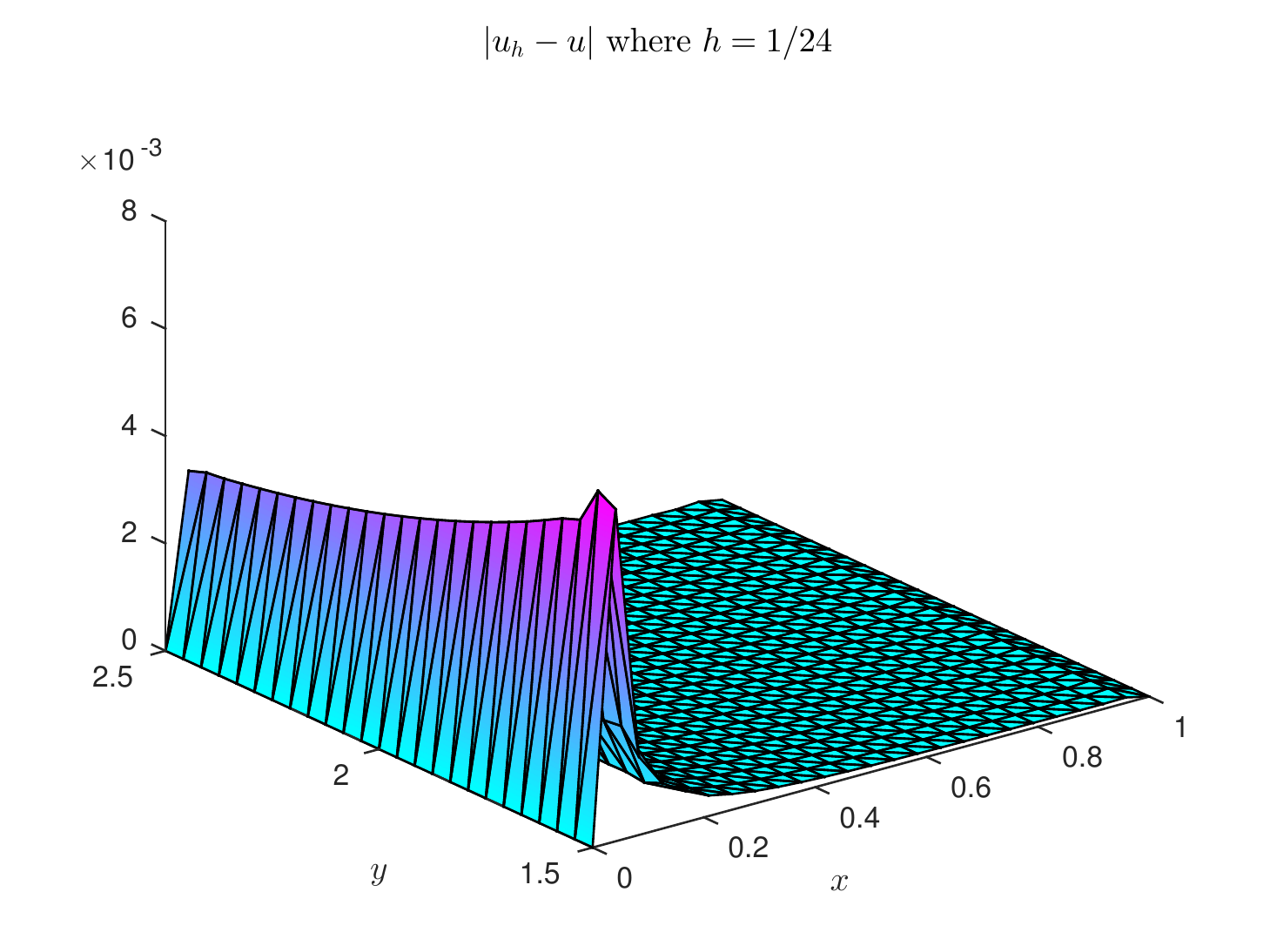}
\caption{The graphs of the error function $|u-\tilde{u}_h|$ (left column) and the error function 
$|u-u_h|$ (right column) with $\a=\frac16$ for $h=\frac16,\frac{1}{12}$, and $\frac{1}{24}$. All 
minimizations were done by using the MATLAB minimization function \texttt{fminunc}.}%
\label{fig5:5}%
\end{figure}

Figure \ref{fig5:5} presents the error plots of both $|\hat{u}_h-u|$ and $|u_h-u|$ over the domain $\W$. 
We observe that $|\tilde{u}_h-u|$ does not converge to zero while $|u_h-u|$ does. In addition, 
Table \ref{table5:3} shows that the cut-off procedure is sufficient in order to guarantee convergence 
for \eqref{eqn5:3}.  Using the same reasoning as to show \eqref{eqn2:7} for Mani\'a's problem, 
a value of $\a<\frac{3}{14}$ is sufficient for the proposed enhanced finite element method to work.  However, 
computing the functional values $\mJha(I_h u)$ with different values of $\a$ shows that 
$\a=\frac12$ and $\a=\frac{10}{17}$ also result in convergent methods.  
\begin{table}[h]
\begin{center}
\begin{tabular}{| c | c | c | c | }\hline
$ h $ & 1/6 & 1/12 & 1/24   \\ \hline  
$\mJ(\tilde{u}_h)$   & 14.84 & 5.71 & 3.21  \\ \hline
$\mJha(\tilde{u}_h)$ & 11.46 & 4.74 & 2.68    \\ \hline
$\mJ(u_h)$  & 3330 & 3914 & 4047    \\ \hline
$\mJha(u_h)$ & 1.28e-1 & 5.45e-3 & 5.28e-4  \\ \hline
\end{tabular}
\caption{The functional values $\mJ$ and $\mJha$ at $\tilde{u}_h$ and $u_h$, where $\tilde{u}_h$ 
and $u_h$ satisfy \eqref{eqn5:4} and \eqref{eqn5:5} respectively for problem \eqref{eqn5:3}. 
Here $\a=\frac16$.}
\label{table5:3}
\end{center}
\end{table}

\section{Conclusion}\label{section5}
In this paper we proposed an enhanced finite element method for variational problems that exhibits 
the {\em Lavrentiev gap phenomenon}. The method uses the Lagrange finite element spaces 
and a discrete energy functional that is constructed by using a novel cut-off procedure. It is simple to 
construct and easy to implemented simply by adding a few extra lines of code. Unlike some existing numerical methods, 
the formulation of the proposed enhanced finite element method does not depend on {\em a priori} knowledge 
about the exact solutions/minimizers. Only one problem-dependent parameter needs to be chosen in order to use 
the method. A sufficient and easy-to-use condition was provided for its determination. We presented 
extensive numerical experiment results for two benchmark problems for the {\em Lavrentiev gap phenomenon}, 
namely, Mani\'a's 1-D problem and Foss' 2-D problem. Our numerical results show that the 
proposed enhanced finite element method works very well, it is effective, robust and convergent. 

It is clear that the crux of the proposed enhanced finite element method is the cut-off procedure, 
it caps the growth of $\nabla v_h$ to the order of $O(h^{-\alpha})$ in the whole domain. Clearly, 
this procedure is independent of the underlying finite element method, consequently, it is natural 
and temptable to substitute the finite element space $V^h_r$ in the formulation by other approximation 
spaces, such as nonconforming finite element spaces and discontinuous Galerkin (DG) spaces. 
After taking care a few details which are well-known when transitioning from the finite element 
method to DG and nonconforming methods, this generalization indeed can be easily done.
In fact, although we only presented the formulation and numerical experiments in the context of 
the finite element method, we have also numerically tested DG methods with and without the cut-off 
procedure, the outcome is exactly same as for the finite element method, that is, the standard DG method
does not give a convergent method, but the enhanced DG method (with the cut-off procedure) does. 
Further details on various generalizations of the work presented in this paper will be reported in
a forthcoming paper. 
  
Besides the generalizations of using approximation spaces other than finite element spaces, 
there are a few important and interesting issues which need to be addressed. The foremost 
issue perhaps is to provide a qualitative convergence analysis for the proposed enhanced finite 
element method. Such a project has already been undertaken using the $\Gamma$-convergence approach
and will be reported later in \cite{FS2016}. 
Since the singularities of the minimizers are often isolated, it is expected and also verified by 
our numerical experiments that the biggest errors are occurred near the singularity points, and 
very fine meshes are required to resolve these singularities. To improve efficiency and to reduce the 
computational cost, It is natural to incorporate unstructured meshes and adaptive algorithms, which
only use very fine meshes near the singularity points, into the enhanced finite element method.
In this regard, the cut-off procedure provides an immediate {\em a posteriori} indicator 
for mesh refinement, that is, the mesh is only refined where the cut-off function is triggered. 
Such an idea is worthy of further investigation and will be reported in another forthcoming work. 



\begin{thebibliography}{00}
\bibliographystyle{abbrv}

\bibitem{BL2006}
Y.~Bai and Z.-P. Li.
\newblock A truncation method for detecting singular minimizers involving the
  {L}avrentiev phenomenon.
\newblock {\em Math. Models Methods Appl. Sci.}, 16(6):847--867, 2006.

\bibitem{BK1987}
J.~M. Ball and G.~Knowles.
\newblock A numerical method for detecting singular minimizers.
\newblock {\em Numer. Math.}, 51(2):181--197, 1987.

\bibitem{CO2010}
C. Carstensen and C. Ortner. 
\newblock Analysis of a class of penalty methods for computing singular minimizers. 
\newblock {\em Comput. Methods Appl. Math.}, 10(2):137-163, 2010.

\bibitem{CO2009}
C. Carstensen and C. Ortner. 
\newblock Computation of the Lavrentiev phenomenon. 
\newblock {\em OxMoS Preprint}, No. 17,  2009.

\bibitem{FS2016}
X. Feng and S. Schnake.
\newblock $\Gamma$-convergence of an enhanced finite element method for approximating 
singular minimizers. 
\newblock {\em in preparation}.

\bibitem{Foss2003}
M.~Foss.
\newblock Examples of the {L}avrentiev phenomenon with continuous {S}obolev exponent dependence.
\newblock {\em J. Convex Anal.}, 10(2):445--464, 2003.

\bibitem{FHM2003}
M.~Foss, W.~J. Hrusa, and V.~J. Mizel.
\newblock The {L}avrentiev gap phenomenon in nonlinear elasticity.
\newblock {\em Arch. Ration. Mech. Anal.}, 167(4):337--365, 2003.

\bibitem{Li1995}
Z.-P. Li.
\newblock A numerical method for computing singular minimizers.
\newblock {\em Numer. Math.}, 71:317--330, 1995.

\bibitem{Li1992}
Z.-P. Li.
\newblock Element removal method for singular minimizers in variational
  problems involving {L}avrentiev phenomenon.
\newblock {\em Proc. Roy. Soc. London Ser. A}, 439(1905):131--137, 1992.

\bibitem{Lavrentiev1927}
M. Lavrentiev.
\newblock Sur quelques problemes du calcul des variations..
\newblock {\em Ann. Mat. Pura e App.}, 4:7--28, 1926.

\bibitem{Mania1934}
B.~Mani\'{a}.
\newblock Soppa un esempio di lavrentieff.
\newblock {\em Ball. Unione Mat. Ital}, 13:147--153, 1934.

\bibitem{Ortner2011}
C.~Ortner.
\newblock Nonconforming finite-element discretization of convex variational
  problems.
\newblock {\em IMA J. Numer. Anal.}, 31(3):847--864, 2011.

\bibitem{Winter1996}
M.~Winter.
\newblock Lavrentiev phenomenon in microstructure theory.
\newblock {\em Electron. J. Diff. Eqns.}, 6:1--12, 1996.

\end{thebibliography}
\end{document}